\documentclass[11pt]{article}
\usepackage[left=1in, right=1in, top=1in]{geometry}
\usepackage[english]{babel}
\usepackage[utf8]{inputenc}

\usepackage{amsmath}
\usepackage{amssymb}
\usepackage{amsfonts}
\usepackage{amsthm}
\usepackage{mathtools}
\usepackage{mathrsfs}

\usepackage[colorlinks, citecolor=blue]{hyperref} 
\usepackage{bbold}

\newcommand{\omprime}{\omega_{\mcx/M}'}
 
\newcommand{\frm}{\mathfrak{M}}

\newcommand{\mc}[1]{\mathcal{#1}}
\newcommand{\on}[1]{\operatorname{#1}}
\newcommand{\pinv}{\pi^{-1}}

\newcommand{\ovm}{\ov{M}}

\newcommand{\rr}{\mathbb{R}}

\newcommand{\ov}[1]{\overline{#1}}

\newcommand{\z}{\mathbb{Z}}
\newcommand{\bos}{\on{bos}}
\newcommand{\ot}{\mc{O}_T} 

\newcommand{\bb}{\mathbb{C}}
\newcommand{\Spec}{\on{Spec}}

\newcommand{\Ber}{\operatorname{Ber}}

\newcommand{\lgr}{\overset{\sim}{\longrightarrow}}

\newcommand{\mfr}{\mathfrak{M}_{g,0,2r}}

\newcommand{\pist}{\pi_*}

\newcommand{\mcx}{\mc{X}} 
\newcommand{\ox}{\mc{O}_X}

\newcommand{\pons}{p_1^*}
\newcommand{\pows}{p_2^*}

\newcommand{\fermp}{P_-}

\newcommand{\partz}{\frac{\partial}{\partial z} } 
\newcommand{\parth}{\frac{\partial}{\partial \theta}}

\newcommand{\dd}{\mc{D}}
\newcommand{\bracket}{ [ \ , \ ]}
\newcommand{\mcz}{\mc{Z}}

\newcommand{\urr}{U_{\rr}}

\newcommand{\lra}{\longrightarrow}

\newcommand{\tbold}{\mc{T}}

\newcommand{\xra}{\xrightarrow}
\newcommand{\deltpull}{\Delta^*}
\newcommand{\vertu}{\vert_U}

\newcommand{\mrr}{M_{\rr}}

\newtheorem{theorem}{Theorem}[section]

\newtheorem{lemma}[theorem]{Lemma}

\newtheorem{proposition}[theorem]{Proposition}
\newtheorem{corollary}[theorem]{Corollary}

\newtheorem{remark}[theorem]{Remark}

\usepackage{appendix} 

\title{ A measure on the moduli space of super Riemann surfaces \\ with Ramond punctures} 
\author{Ron Donagi, Nadia Ott}
\date{{}}

\begin{document}

\maketitle

\section{Introduction}

In perturbative string theory, an essential ingredient is a certain measure on the ordinary moduli space $\mc{M}_g$ of curves. At a point corresponding to a curve $C$, the Mumford isomorphism relates the determinant of the tangent space of $\mc{M}_g$ to a power of the determinant of cohomology of the canonical bundle $\omega_C$ on the curve. The cohomology of the canonical bundle $\omega_C$, and therefore also the determinant of this cohomology, has a natural Hermitian metric given by integration on $C$, and via the Mumford isomorphism this determines a measure on $\mc{M}_g$. We can think of this Hermitian metric as a pairing between a cohomology on $C$ and the corresponding cohomology on the complex-conjugate curve $\overline{C}$.

An equivalent way of expressing this is in terms of the period map. This sends $H^0(C,\omega_C)$ to $H^1(C,\bb)$. So the standard intersection form on $H^1(C,\bb)$ pulls back to a Hermitian metric on $H^0(C,\omega_C)$ and therefore also on its determinant.

 In superstring theory, this measure is generalized to the  \emph{supermeasure}, a measure on the moduli space $\mathfrak{M}_g$ of genus $g$ super Riemann surfaces. The construction uses super versions of the previous ingredients: a super Mumford isomorphism, reviewed below \ref{SRS}, and a Hermitian metric on the determinant of cohomology on the supercurve, which again can be interpreted in terms of a super period map. This provides a holomorphic measure on an open subset (see below) of the product of 
$\mathfrak{M}_g$ with an appropriate complex-conjugate version $\overline{\mathfrak{M}_g}$.
In a final step, this is restricted to the integration cycle needed in superstring theory. We will say a few words about this below.

Unlike the classical period matrix, the normalized super period matrix is in general only meromorphic along the bad locus. The unnormalized period morphism used in this paper is holomorphic as a morphism of sheaves, but along the bad locus the relevant cohomology jumps and the period image is no longer a vector bundle of the generic rank. Witten conjectured that the supermeasure extends smoothly across this locus for genus $\leq 11$. This was rigorously proved in \cite{felder2019regularity}. Very recently, Deligne extended this result to all genera, as reported in 
\cite{felder2024superperiods}.

For further applications to perturbative superstring theory,  measures on the moduli spaces of super Riemann surfaces with punctures are needed. There are two types of such punctures: Neveu-Schwarz (NS) and Ramond.  The case of Neveu-Schwarz punctures, described in \cite{smiwithRamondwitten2015notes} and Appendix \ref{appendix: ns punctures} of this paper, is relatively straightforward since the moduli space of super Riemann surfaces  with Neveu-Schwarz punctures admits a projection onto the moduli space of unpunctured super Riemann surfaces. 

The situation is very different though for Ramond punctures: Ramond punctures are divisors along which the superconformal structure  degenerates, and there is no sense in which they can be "forgotten" to produce a projection onto the unpunctured supermoduli space. 

In this work we define a measure on the moduli space $\mfr$ of super Riemann surfaces with $2r$ Ramond punctures (the number of such punctures is always even), generalizing the supermeasure on the unpunctured supermoduli space. The definition of the measure on $\mfr$ uses the adaptation to the Ramond case of the familiar ingredients: the Mumford isomorphism, described in Appendix C of \cite{smiwithRamondwitten2015notes},
and a generalization of the super period map to $\mfr$, defined in Section \ref{sectiontitle: super period map} of this paper.  As in the unpunctured case, we define this measure initially away from a certain bad locus. We are able to prove that, for $r \geq 2$, it extends smoothly across the bad locus.  The basic reason that the measure extends across the bad locus is Hartogs-like: for $r \geq 2$ we show that the bad locus has codimension $\geq 2|0$. The case $r=0$ is Deligne’s result, and the case $r=1$ remains open.

The measure we construct is defined a priori on an open
neighborhood of a Zariski-open subset of the diagonal. Our main result, Theorem \ref{main}, is that for $g\ge2$ and $r>1$ it extends to an open neighborhood of the entire diagonal. The extended measure is allowed to have zeros; indeed, by Theorem \ref{exactcodim} together with Proposition \ref{not volume form}, the extended measure necessarily has zeros along the bad locus whenever $r>2$ and the bad locus is nonempty.
 
In the product $\mathfrak{M}_{g,0,2r} \times \overline{\mathfrak{M}_{g,0,2r}}$, physicists consider the quasidiagonal - roughly, the locus of pairs $X,\widetilde{X}$ where the bosonic curves underlying $X,\widetilde{X}$ are complex conjugates of each other. This contains the diagonal, but is strictly bigger. The integration cycle (for type II superstring theory) is a thickening of this quasidiagonal. It is not clear whether our measure extends to an open neighborhood of the quasidiagonal.

In the remainder of this Introduction we outline our argument in more detail.

\subsection{Super Riemann surfaces and their moduli}\label{SRS}
We will work mostly with super Riemann surfaces with Ramond punctures. Their definitions and elementary properties are reviewed in Section \ref{sectiontitle: Ramond Punctures}. Throughout the body of the paper we assume $g\ge2$. The moduli of smooth super Riemann surfaces of genus $g$ with $2r$ Ramond punctures is a smooth Deligne--Mumford superstack $\frm = \frm_{g,0,2r}$ \cite{felder2020moduli}. (The genus-zero case is treated in \cite{ott2023supermoduli}, but is not needed below.) For simplicity, we refer to this superstack as the moduli space. A morphism $T \to \frm$ from a superscheme $T$ determines a family over $T$ of genus $g$ super Riemann surfaces with Ramond divisor $R$ consisting of $2r$ Ramond punctures, each supported on an irreducible component $R_i$. There is a universal family $\pi: \mcx \to \frm$ with the property that, for every morphism $T \to \frm$, the pullback $\pi_T: \mcx \times_{\frm} T \to T$ is the corresponding family of super Riemann surfaces with Ramond punctures over $T$. By a point of $\frm$ we mean a morphism $\Spec \bb \to \frm$.

Near each irreducible component $R_i$ of the Ramond divisor, we can choose local coordinates $(z,\theta)$ on $X$ such that $R$ is locally defined by $z = 0$, and the distribution $\mc{D}$, which is maximally non-integrable away from $R$, is generated by the odd vector field $D_{\theta} := \parth + z\theta \partz$. Such coordinates are called \emph{superconformal}. For a family of super Riemann surfaces with Ramond punctures $X \to T$, superconformal coordinates exist \'etale locally on $X$ \cite[Proposition 3.6]{bruzzo2021supermoduli}.

If $(z', \theta')$ is another set of superconformal coordinates near a Ramond puncture, then
\begin{equation} \label{orientation}
\theta' = \pm(\theta + \tau) \mod z,
\end{equation}
for some odd section $ \tau $ on the base superscheme $T$, cf. \cite{witten2015super}. A reduction of this ambiguity to just translations by $\tau$ is called an {\em{orientation}}. We let 
$\widetilde{\frm} = \widetilde{\frm_{g,0,2r}}$ 
denote the moduli space of  super Riemann surfaces of genus $g$ with $2r$ {\em oriented} Ramond punctures. It is a covering of $\frm$ of degree $2^{2r}$.

In addition to the Berezinian, or dualizing sheaf, $\omega_X$, we will use the extended Berezinian sheaf $\omega'_X$ defined in \cite{witten2015super}. By definition, sections of $\omega'_X$ are sections of $\omega_X(R)$ whose residue along the Ramond divisor is constant. Both $\omega_X$ and $\omega'_X$ are canonically defined for every family of super Riemann surfaces with Ramond punctures. In particular, they are defined on the pullback of the universal curve to any \'etale cover of $\mcx$, and these local constructions agree on overlaps. Hence they descend to sheaves on $\mcx$, and therefore pull back to sheaves on any family of super Riemann surfaces with Ramond punctures. We define the relative version of the extended Berezinian sheaf in Section \ref{sectiontitle: de Rham with Rmond}.

A family of super Riemann surfaces of genus $g$ with $2r$ Ramond punctures over a purely bosonic base scheme (e.g., $T=\Spec \bb$) can be described by the bosonic data $(C,D,L)$ where $C$ is a smooth, compact Riemann surface of genus $g$, $D=p_1 + \dots + p_{2r}$ an effective divisor on $C$, and $L$ is a line bundle equipped with an isomorphism $i:L^2\cong\omega_C(D)$. We refer to $C,D$ and $L$ as the underlying bosonic curve, bosonic Ramond divisor, and twisted spin structure, respectively. An orientation on $X$ amounts to choosing a trivialization of $L \vert_D$ whose square equals the composition 
$\operatorname{Res}_D \circ i: L^2 \to \bb$ 
of the given isomorphism $i$ with the residue along the divisor $D$.

In this case, we have the following useful decompositions of $\omega_X$ and $\omega_X'$ in terms of purely bosonic data: 
\footnote{Let $i:C \to X$ denote the inclusion. The sheaf on the right is a sheaf of $\mc{O}_C$-modules, so the object on the left should really be $i^{-1}{\omega_X}$. We ignore this very slight inaccuracy here and elsewhere in this paper.}
\[
{\omega_X} = \omega_C \oplus \Pi L(-D), \ \ \ \ {\omega'_X} = \omega_C \oplus \Pi L.
\]

\noindent From this decomposition we see that, the dimensions of the cohomology groups are of the form:
\begin{equation}\label{dims}
h^0(\omega_X)=g|\epsilon, \ h^1(\omega_X)=1|r+\epsilon, \ 
h^0(\omega'_X)=g|r+\epsilon, \ h^1(\omega'_X)=1|\epsilon,
\end{equation}
with $\epsilon:= h^0(L)-r=h^0(\omega_C \otimes L^{-1}) =h^0(L(-D))$ and $\epsilon=0$ for generic choices. The {\em bad locus} mentioned above is defined as the closed subspace of $\frm$ where $\epsilon >0$.
For more details, see \cite{witten2015super, donagi2023bad}

\subsection{Mumford isomorphisms}
The original Mumford isomorphism (Theorem 5.10 in \cite{mumford1977stability}) is
\[
L_2 \cong {L_1}^{13},
\]
where  $\pi: \mc{C} \to \mc{M}$ denotes the universal curve over the moduli space of curves of genus $\geq 2$, $\omega := \omega_{\mc{C} / \mc{M}}$, and
$L_i := \det R\pi_* (\omega^{\otimes i})$. It is often used together with the  identification of $R\pi_* (\omega^{\otimes 2})$ with the cotangent bundle of $\mc{M}$, hence of $L_2$ with the canonical bundle of $\mc{M}$.

The super Mumford isomorphism \cite{voronov1988formula, rosly1989superconformal} is

\begin{equation}\label{MI}
\Ber_3 \cong {\Ber_1}^5,
\end{equation}
where now $\pi: \mcx \to \frm$ denotes the universal supercurve over supermoduli, 
$\omega := \omega_{\mcx / \frm}$ is the (relative) dualizing sheaf, and
\begin{equation}\label{Ber}
    {\Ber}_i := \Ber R\pi_* (\omega^{\otimes i}).
\end{equation}
This is often used together with the  identification of $ R\pi_* (\omega^{\otimes 3})$ with the cotangent bundle of $\frm$, hence of 
$Ber_3$ with the dualizing sheaf of $\frm$.


On $\frm_{g,0,2r}$, i.e. in the presence of Ramond punctures, the
cotangent theory involves the twist by twice the Ramond divisor. More
precisely, Serre duality identifies the canonical line of supermoduli with
\[
\omega_{\frm_{g,0,2r}}
\cong
\Ber\!\left(R\pi_*\omega_{\mc X/\frm}^{3}(2R)\right).
\]
The twist by $2R$ does not change the Berezinian of cohomology. Indeed,
there is an exact sequence
\[
0\longrightarrow
\omega_{\mc X/\frm}^{3}
\longrightarrow
\omega_{\mc X/\frm}^{3}(2R)
\longrightarrow
\mc J_R
\longrightarrow0,
\]
where $\mc J_R$ is the sheaf of polar parts along the Ramond divisor.
Appendix C of \cite{smiwithRamondwitten2015notes} gives a canonical
trivialization
\[
\Ber(R\pi_*\mc J_R)\cong\mc O_{\frm}.
\]
Consequently
\[
\Ber\!\left(R\pi_*\omega_{\mc X/\frm}^{3}(2R)\right)
\cong
\Ber\!\left(R\pi_*\omega_{\mc X/\frm}^{3}\right).
\]
Combining this with the super Mumford isomorphism gives
\begin{equation}\label{BerM}
\omega_{\frm_{g,0,2r}}
\cong
\Ber^5\!\left(R\pi_*\omega_{\mc X/\frm}\right).
\end{equation}

\subsection{The pairing}

The symmetry hinted at by equation \eqref{dims} is not accidental. 
The Berezinians of cohomology of $\omega$ and $\omega'$ are actually isomorphic, cf. Theorem \ref{firstresult}: 
What is changing in the cohomology as we move from $\omega$ to $\omega'$ is that an $(r+\epsilon)$-dimensional piece disappears from the odd part of $H^1$, and is replaced by its dual vector space that appears in the odd part of $H^0$. This has no effect on the Berezinian.

So instead of constructing a pairing over supermoduli involving 
${\Ber} := \Ber R\pi_* (\omega),$
as in \eqref{Ber}, we can 
construct a pairing involving the isomorphic
${\Ber'} := \Ber R\pi_* (\omega').$  So our task is to construct a natural pairing on $\Ber R\pi_* (\omega').$

We start by defining some local systems on supermoduli space $\frm:=\frm_{g,0,2r}$, where we will use $\pi: \mcx \to \frm$ to denote the universal family with Ramond divisor $R$. One is 
\[
\Lambda_0 := R^1\pi_* \z.
\]
This has rank $2g|0$ and carries the standard $\z$-valued symplectic pairing corresponding to the intersection pairing on the underlying curve $C$. The pairing is invariant under the modular group.

Another local system comes from the sheaf $\mc{Z}_{R/\frm}^1$ of relative closed-one forms on $R$. Let $i: R \to \mcx$ denote the natural inclusion.  There is a natural local system $\Lambda_1$ of free abelian groups of rank $0|2r$ on $\frm$, carrying a non-degenerate symmetric bilinear pairing $J_1: \Lambda_1 \times \Lambda_1 \to \z$. This pairing is invariant under the group $(\z/2\z)^{2r}$ of changes of \emph{orientation on $R$}, cf. Section \ref{sectiontitle: de Rham with Rmond}. Since $R$ is smooth and finite of degree $2r$ over $\frm$,  $\pist i_* \mc{Z}_{R/\frm}^1$ is locally free of rank $0|2r$, and we show in Section \ref{section with orientation}, 
that $\pist \mc{Z}_R^1 \cong \Lambda_1 \otimes_{\z} \mc{O}_{\frm_{g,0,2r}}$.

We set 
\[ \Lambda:= \Lambda_0 \oplus \Lambda_1, \qquad \Lambda_{\frm}:= \Lambda \otimes \mc{O}_{\frm}\] 
This is a local system of free abelian groups of rank $2g|2r$ with a non-degenerate pairing, 
\[
J:=J_0 \oplus J_1: \Lambda \times \Lambda \to \z, \qquad J_{\frm}:= J \otimes \mc{O}_{\frm}
\]
where $J_0$ is the alternating intersection pairing on the even part $\Lambda_0$, and $J_1$ is the symmetric pairing on the odd part $\Lambda_1$.



We define the super period map $P_{X/T}$ for families $X \to T$ of super Riemann surfaces with Ramond punctures in  Section \ref{sectiontitle: de Rham with Rmond}. For the universal family $\pi : \mcx \to \frm$, the super period map is
\begin{equation}\label{periodmap}
P: \pi_* \omega_{\mcx/\frm}' \to \Lambda \otimes \mc{O}_{\frm}.
\end{equation}
The map $P_{\mcx/\frm}$ is the direct sum of two maps
\[ P^+: \pi_* \omega_{\mcx/\frm}' \to \Lambda_0 \otimes \mc{O}_{{\frm}}
\]  
\[
P^-:\pi_* \omega_{\mcx/\frm}' \to \Lambda_1 \otimes \mc{O}_{\frm}
\]
which we call the even and odd period maps, respectively. Both maps are grading-preserving.

The $2g$ even periods of a section of $\pist \omega'$ are its (Berezin)-integrals over a basis of $1$-cycles in cohomology, while its $2r$ odd periods are its residues along the $2r$ Ramond punctures. 

We show in Theorem \ref{bigone} that for any super Riemann surface $X$ over $\Spec \bb$ such that $\Spec \bb \to \frm$ is in the good locus, its super period map $P_X: H^0(X, \omega_X') \to \Lambda_{\bb}$ is an injective morphism of vector bundles whose image $\on{im}(P_X)$ is $J$-Lagrangian. 
Furthermore, since $\Lambda_{\bb}$ is the complexification of the lattice $\Lambda$, we can define a complex conjugate map $\lambda \otimes v \to \lambda \otimes \ov{v}$, and thus also take the complex conjugate of subspaces in $\Lambda_{\bb}$. In Theorem \ref{bigone}, we show the restriction of $J$ to $\on{im}(P) \times \ov{\on{im}(P_X)}$ is non-degenerate.

\subsection{The measure}

Throughout this section we set $\frm:=\mfr$. 
 We want to define the pairing: 
\begin{equation} 
\label{bilinearpairing} \pons \pist \omega' \otimes \pows \ov{\pist \omega'} \overset{\pons P \otimes \pows \ov{P}}
{\xhookrightarrow{\hspace{1cm}}} p_1^{-1}\Lambda \otimes p_2^{-1} \Lambda \otimes  \mc{O}_{\frm \times \ov{\frm}} \overset{J \otimes 1}{\longrightarrow} \mc{O}_{\frm \times \ov{\frm}}.
\end{equation}

\noindent In order for the second map to make sense, we need to be able to identify $p_1^{-1}\Lambda_1$ and $p_2^{-1} \Lambda_1$. We do this by restricting everything in \eqref{bilinearpairing} to an open neighborhood $\frm'$ of the diagonal $\Delta(\frm) \subset \frm \times \ov{\frm}$.

We want to take the Berezinian of the resulting  pairing. In order to do that, we need the pairing to be non-degenerate. We first prove this away from the bad locus (cf. Corollary \ref{corollary of big one}). Or more precisely, away from the bad diagonal, $\Delta(\mc{B}) \subset \Delta(\frm)$.   
Combining the pairing on the Berezinian of \eqref{bilinearpairing} with the super Mumford isomorphism, $\omega_{\frm}= \Ber^5(R \pist \omega_{\mcx/\frm})$, and with the identification of $\Ber^5(R \pist \omega')$ with $\Ber^5(R \pist \omega)$, gives the measure on the complement of the bad diagonal.

In Theorem \ref{dimensioncounts} we show that the bad locus has complex codimension at least $2$. In Section \ref{sectiontitle: extension across the bad locus} we do not attempt to extend the period image as a vector bundle. Instead, we extend directly the inverse Berezinian section which defines the measure. The extension is first obtained across the degeneracy locus of the period pairing by a local block-matrix formula for the Berezinian, and is then extended across the two factors of the bad locus by Hartogs' theorem. This proves that the measure extends smoothly to an open neighborhood containing the entire diagonal. The extension need not be nowhere vanishing.

\subsection{Further comments}

What we call Ramond punctures corresponds, in the superstring literature, to Ramond-Ramond punctures, indicating that both $\frm$ and $\ov{\frm}$ parametrize super Riemann surfaces with Ramond punctures. Other possibilities are NS-NS punctures and the mixed NS-R and R-NS punctures. We use the simplified terminology since R-R punctures are the only ones we consider here.

It is not clear how or if the measure constructed in this paper relates to the one needed in perturbative superstring theory. 
For instance, in the unpunctured case, the correct  integration cycle  is not the diagonal, but rather a  larger subset $\Gamma$ of $\frm \times \ov{\frm}$ called the quasi-diagonal. The quasi-diagonal contains the diagonal, but it is strictly bigger.

In the Ramond case, the quasi-diagonal $\Gamma \subset \mfr \times \ov{\mfr}$
parameterizes those pairs of super Riemann surfaces with Ramond punctures whose bosonic data $(C,D,L),(C',D',L')$
are such that $C'$ is the complex conjugate of $C$, $D'$ is the complex conjugate of $D$, but
no constraint is imposed on the twisted spin structure $L'$.
Over the quasi-diagonal, and hence over a small neighborhood of the quasi-diagonal, the even lattices $\Lambda_0, \Lambda'_0$ are identified, so the pairing \eqref{bilinearpairing} is well-defined. To identify the lattices $\Lambda_1$ and $\Lambda_1'$, we would need to pass to the cover
$\widetilde{\frm} \times \widetilde{\overline{\frm}}$,
parameterizing pairs   $(X, \ov{X})$ of super Riemann surfaces of genus $g$ with $2r$ {\em oriented} Ramond punctures.  It is a covering of $\frm \times \ov{\frm} $ of degree $2^{4r}$. Dividing by the simultaneous action of changes-of-orientation on the two factors, we get an intermediate cover 
$\widetilde{\frm \times \ov{\frm}}$, of degree $2^{2r}$.

The pull-back of the super period maps $P$ and $\ov{P}$ to $\widetilde{\frm}$ and $\widetilde{\ov{\frm}}$ are essentially the same as the versions on $\frm$ and $\ov{\frm}$, except that
 $\Lambda_1$ and $\Lambda_1'$ can now both be identified  with the trivial local system $V:= \underline{\z}^{0|2r}$,  and $P_o$ is composed with the trivialization over the Ramond divisor specified by the orientation:
\[
\widetilde{P_o}:\widetilde{\pi_* \omega'}  \to V \otimes \mc{O}_{\widetilde{\frm_{g,0,2r}}}.
\] 

Pulling back all terms in \eqref{bilinearpairing} to $\widetilde{\frm} \times \widetilde{\overline{\frm}}$ or to $\widetilde{\frm \times \ov{\frm}}$, and composing with the identification of $\Lambda_1$ and $\Lambda_1'$ with $V$, we get a well-defined pairing: 
\begin{equation} 
\label{coverbilinearpairing} \widetilde{\pons \pist \omega'} \otimes \widetilde{\pows \ov{\pist \omega'}} \overset{ \widetilde{\pons P} \otimes \widetilde{\pows \ov{P}}}
{\xhookrightarrow{\hspace{1cm}}} p_1^{-1}\Lambda \otimes p_2^{-1} \Lambda \otimes  \mc{O}_{{\widetilde{\frm \times \ov{\frm}}}} \overset{J \otimes 1}{\longrightarrow} \mc{O}_{\widetilde{\frm \times \ov{\frm}}}.
\end{equation}

Various issues remain unclear. Is the restriction of this pairing to the quasi-diagonal $\widetilde{\Gamma}$ in $\widetilde{\frm \times \ov{\frm}}$ non-degenerate? Does it descend to a neighborhood of the quasi-diagonal in the original
$\mfr \times \ov{\mfr}$? Is there a natural way to introduce a dependence on the spinor indices of Ramond vertex operators? We hope to return to these questions elsewhere.

\subsection*{Acknowledgements}
We thank Jacques Distler, David Kazhdan, Sasha Polishchuk, Eric Sharpe, Vivek Shende, Sasha Voronov,   and especially Edward Witten for very helpful discussions. During the preparation of this work, RD's research was supported in part by NSF grants DMS 2001673 and 2401422, by NSF FRG grant DMS 2244978, and by Simons HMS Collaboration grant 390287.

\section{Preliminaries: Super Riemann surfaces with Ramond punctures}
\label{sectiontitle: Ramond Punctures}

Fix $g \ge 2$.
\smallskip

Throughout this paper, we use \emph{supercurve} to mean a smooth, proper, connected superscheme $X$ over $\Spec \bb$ of dimension $1|1$ whose bosonic reduction $C:=X_{\bos}$ is an ordinary curve. The genus $g$ of $X$ is the genus of $C$. A family of genus $g$ supercurves over a superscheme $T$ is a smooth, proper morphism $\pi: X \to T$ of relative dimension $1|1$ such that each fiber $X_t$ is a genus $g$ supercurve over $\Spec k(t)$.

We will write $X$ for any family of supercurves over $T=\Spec \bb$. For any point $t \in T$, we write $k(t)$ for the residue field and denote the fiber by $X_t$.

We recall that a family of genus $g$ \emph{super Riemann surfaces} is a family of genus $g$ supercurves $X \to T$ equipped with a \emph{superconformal structure}, that is, a rank $(0|1)$ subbundle $\mc{D} \subset \mc{T}_{X/T}$ of the relative tangent sheaf which is maximally non-integrable in the sense that the composition
\begin{equation} \label{beta comp}
\mc{D}^{\otimes 2} \xrightarrow{[\ , \ ]} \tbold_{X/T} \to \tbold_{X/T} / \mc{D},
\end{equation}
where $\bracket$ is the supercommutator, is an isomorphism of locally free sheaves on $X$. We then have the following short exact sequence of locally free $\ox$-modules:
\begin{equation}
    \label{sesofdist}
    0 \lra{{}} \mc{D} \lra{{}} \tbold_{X/T} \lra{{}} \dd^{\otimes 2} \lra{{}} 0.
\end{equation}

\begin{remark}
The above definition of a super Riemann surface is from the algebro-geometric perspective, in which a supercurve is a superscheme of dimension $1|1$. Since the individual supercurves are proper, super GAGA allows us to pass between their algebraic and analytic descriptions. For the (non-proper) smooth moduli stack itself, when analytic arguments are used below we work with its analytification rather than invoking GAGA for the moduli space.
\end{remark}

Super Riemann surfaces have two types of punctures: Neveu-Schwarz (NS) punctures and Ramond punctures. Here we will be interested only in Ramond punctures. The case of NS punctures is considered in Appendix \ref{appendix: ns punctures}.

A \emph{family of genus $g$ super Riemann surfaces with $n_R$ Ramond punctures} is a family of genus $g$ supercurves $X \to T$ equipped with the following additional data:
\begin{enumerate}
    \item A relative effective Cartier divisor $R \subset X$, smooth over $T$, whose bosonic reduction $R_{\bos} \subset C=X_{\bos}$ is finite \'etale of degree $n_R$ over $T_{\bos}$. The divisor $R$ is called the \emph{Ramond divisor}. We assume its components are labeled $R_1, \dots, R_{n_R}$ and refer to these as the \emph{Ramond punctures}.
    \item A rank $(0|1)$ subbundle $\mc{D} \subset \tbold_{X/T}$ such that $\dd^{\otimes 2} \xrightarrow{\bracket} \tbold_{X/T} \to (\tbold_{X/T}/\dd)(-R)$ is an isomorphism of locally free sheaves on $X$. In this case, the superconformal structure $\mc{D}$ fits into the following short exact sequence of sheaves on $X$:
\begin{equation} \label{unpunctureddefiningsequence}
0 \longrightarrow \mc{D} \longrightarrow \mc{T}_{X/T} \longrightarrow \mc{D}^{\otimes 2}(R) \to 0.
\end{equation}
\end{enumerate}

An important observation is that $\dd \vert_R$ is naturally identified with the relative tangent bundle $\tbold_{R/T}$ on $R$.

An isomorphism of super Riemann surfaces with $2r$ Ramond punctures from $(X/T,R,\mc{D})$ to $(X'/T,R',\mc{D}')$, called a \emph{superconformal isomorphism}, is an isomorphism of supercurves $\phi:X \cong X'$ over $T$ such that $\phi(R_i)=R_i'$ for all $1 \le i \le 2r$, and
\[
d \phi(\mc{D})=\phi^* \mc{D}',
\]
where $d \phi: \mc{T}_{X/T} \lgr \phi^* \mc{T}_{X'/T}$ is the induced isomorphism.

The number $n_R$ of Ramond punctures is always even, and we will henceforth replace $n_R$ by $2r$. This fact is well-known and follows from the correspondence with twisted spin curves. We now recall this correspondence, primarily to set up notation that will be used throughout the paper.

A family of genus $g$ \emph{spin curves with $n_R$ Ramond punctures}, or twisted spin curves for short, is an ordinary family of genus $g$ curves $C \to T$ over an ordinary scheme $T$, together with a relative effective Cartier divisor $D \subset C$ which is finite \'etale of degree $n_R$ over $T$, and a line bundle $L$ with a fixed isomorphism
\[
\alpha:L^{\otimes 2} \cong \omega_{C/T}(D).
\]
Note that $\deg(L)=g-1+n_R/2$, and hence $n_R$ must be even.

Now let $X \to T$ be a family of super Riemann surfaces with $n_R$ Ramond punctures, with $T$ a purely bosonic scheme, and let $C$ denote the underlying bosonic curve. Since $T$ contributes no odd coordinates, $J^{\otimes 2}=0$ and $J$ is a rank $(0|1)$ line bundle on $C$. Restricting the sequence \eqref{unpunctureddefiningsequence} to $C$ gives canonical isomorphisms
\[
J \cong \mc{D}\vert_C^{-1}, \qquad J^{\otimes 2} \cong \mc{D} \vert_C^{-2} \cong \omega_{C/T}(D),
\]
where $D:=R \vert_C$, and all tensor products are over $\mc{O}_C$. Setting $L:= \Pi J$, the data $(C \to T, D=R \vert_C, L)$ defines a twisted spin curve over $T$.

A further fundamental fact about super Riemann surfaces with Ramond punctures is the existence of superconformal coordinates, proved in Lemma 2.6 of \cite{felder2020moduli}. If $p$ is a point of $X$ not contained in the Ramond divisor, then this is the standard statement that there are local coordinates $(z,\theta)$ in which $\mc{D}$ is locally generated by the odd vector field
\[
\parth + \theta \partz.
\]
If $p \in R$, then there are local coordinates $(z,\theta)$ such that $R$ is locally defined by $z=0$ and $\mc{D}$ is generated by
\[
D_{\theta} := \parth + z \theta \partz.
\]
The coordinates $(z,\theta)$ are called \emph{superconformal coordinates}. In the algebro-geometric setting, superconformal coordinates exist \'etale locally on $X$, and in families after passing to an appropriate \'etale double cover.

\subsection{Orientation on the Ramond divisor} \label{section with orientation}


Let $i:R\hookrightarrow X$ be the natural inclusion of the Ramond divisor, and let $(z,\theta)$ be a pair of superconformal coordinates on $X$ such that $R$ is locally defined by $z=0$. Then $i^*(z)=0$, while
\[
x:=i^*(\theta)
\]
is a local odd coordinate on $R$. We refer to a set of global coordinates $(x^1, \dots, x^{2r})$ on $R$ as superconformal coordinates if and only if they are obtained by restriction of superconformal coordinates on $X$, and call an automorphism of $R$ a superconformal automorphism if it sends superconformal coordinates to superconformal coordinates.

In \cite{witten2015super}, Witten shows that if $(x^1, \dots, x^{2r})$ are superconformal coordinates on $R$, and $\phi$ is a superconformal automorphism of $R$, then
\begin{equation} \label{orientationonR}
\phi^*(x^i)=\pm(x^i+\tau_i),
\end{equation}
where $\tau_i$ is an odd function on the base $T$; in particular, $\tau_i=0$ if $T$ is purely bosonic.
\medskip

Let $d_{R/T}$ denote the relative exterior derivative, and set
\[
\mcz_{R/T}^1:=\ker\left(d_{R/T}:\Omega_{R/T}^1 \to \Omega_{R/T}^2\right).
\]
Since $i_*d_{R/T}$ is $\pinv(\ot)$-linear, $i_*\mcz_{R/T}^1$ is a $\pinv(\ot)$-module on $X$. Note that it is \textbf{not} an $\ox$-module.

For each $t \in T$, choose superconformal coordinates $(x^1,\dots,x^{2r})$ on $R_t$. Then the differentials
\[
dx^1,\dots,dx^{2r}
\]
form a basis for $H^0(R_t,\mcz_{R_t}^1)$. Indeed, a global one-form on $R_t$ is locally of the form $v=(a+bx^i)dx^i$ for some $a,b \in k(t)$, and $d(v)=0$ if and only if $b=0$. Thus,
\[
\on{dim}_{k(t)} H^0(R_t,\mcz_{R_t}^1)=0|2r
\]
for all $t \in T$, and $\pist \mcz_{R/T}^1$ is locally free of rank $0|2r$ on $T$. A local frame for $\pist\mcz_{R/T}^1$ can therefore be chosen as the differentials $(dx^1,\dots,dx^{2r})$ of any local set of superconformal coordinates $(x^1,\dots,x^{2r})$ on $R$.

We call any choice of local frame for $\pist \mcz_{R/T}^1$ arising from a choice of superconformal coordinates on $R$ an \emph{orientation} on $R$, and any such local frame a \emph{local orientation} on $R$.

Choose a local orientation $(dx^1,\dots,dx^{2r})$, set
\[
V=\z dx^1+\cdots+\z dx^{2r},
\]
and note that $V \otimes \ot$ is locally isomorphic to $\pist \mcz_{R/T}^1$. From \eqref{orientationonR}, if $(dx'^1,\dots,dx'^{2r})$ is another local orientation, then
\[
dx'^i=\pm dx^i, \qquad \forall \ 1\le i\le 2r.
\]
Hence $V$ and $V'$ agree on overlaps, and therefore glue to a locally constant sheaf $V_{\z}$ on $T$. Thus,
\[
\pist \mcz_{R/T}^1 \cong V_{\z} \otimes \ot.
\]

Define a quadratic form
\[
J_1:V_{\z}\longrightarrow \z
\]
by
\[
J_1\left(\sum_i n_i\,dx^i\right)=\sum_i n_i^2,\qquad n_i\in\z,
\]
equivalently $\langle dx^i,dx^j\rangle=\delta_{ij}$. Since changes of local orientation are given by $dx^i\mapsto \pm dx^i$, the transition functions lie in $O(2r,\z)$, so this quadratic form is independent of the chosen frame and defines a canonical $\z$-valued quadratic form on the local system $V_{\z}$.

\begin{lemma} \label{oddtarget} Let $X$ be a super Riemann surface with  $2r$ Ramond punctures over $\Spec \bb$, let $C$ denote its underlying bosonic curve, $L$ the twisted spin structure, and $D=R \vert_C$. 
Then there is a natural identification
\[
Z_R^1=\Pi L\vert_D, \qquad Z_R^1:=H^0(R,\mcz_R^1).
\]
Furthermore, the complexification of $J_1$ is exactly the residue pairing on $L \vert_D$. For each $p_i \in D$, the residue pairing $L_{p_i} \otimes L_{p_i} \to \bb$ sends $s \otimes s' \in L^2=\omega_C(D)$ to $\on{res}_{p_i}(s \otimes s')$. 
\end{lemma}

\begin{proof}
Since $\ox=\mc{O}_C\oplus\Pi L$, for any $\ox$-module $\mc{M}$ on $X$, we have
\[
\mc{M}=\mc{M}\vert_C\oplus(\mc{M}\vert_C\otimes_{\mc{O}_C}\Pi L).
\]
Applying this to $\dd^{-1}$ and using the canonical isomorphism $\dd^{-1}\vert_C\cong\Pi L$, we find
\[
\dd^{-1}
=
\dd\vert_C^{-1}\oplus(\dd\vert_C^{-1}\otimes\Pi L)
\cong
\omega_C(D)\oplus\Pi L,
\]
where $D=R\vert_C$. Restricting further to $R\subset X$ and using the identification $\mc{D}^{-1}\vert_R=\Omega_R^1$, we obtain
\[
\Omega_R^1=\omega_C(D)\vert_D\oplus\Pi L\vert_D.
\]
Thus,
\[
H^0(R,\Omega_R^1)
=
H^0(D,\omega_C(D)\vert_D)\oplus\Pi H^0(D,L\vert_D).
\]
Since $Z_R^1$ is odd of dimension $0|2r$, we have
\[
Z_R^1\subseteq\Pi H^0(D,L\vert_D).
\]
The right-hand side also has dimension $0|2r$, with one $0|1$ contribution from each of the $2r$ components of $D$, so the inclusion is an equality. Finally, since $D$ is zero-dimensional,
\[
H^0(D,L\vert_D)=L\vert_D,
\]
and therefore
\[
Z_R^1=\Pi L\vert_D,
\]
as claimed.
\end{proof}

\section{Periods of super Riemann surfaces with Ramond punctures}
\label{sectiontitle: de Rham with Rmond}\label{sectiontitle: super period map}

Let $\pi : X \to T$ be a family of genus $g$ super Riemann surfaces with $2r$ Ramond punctures.

We start by defining the relative extended Berezinian sheaf $\omega_{X/T}'$ on the family $\pi: X \to T$. From the identifications $\omega_{X/T} = \dd^{-1}(-R)$ and $\Omega_{R/T}^1 = \dd^{-1} \vert_R$, we get an exact sequence
\[
0 \longrightarrow \omega_{X/T} \longrightarrow \omega_{X/T}(R) \to i_*\Omega_{R/T}^1 \longrightarrow 0.
\]
where $i: R \to X$ is the canonical inclusion of the Ramond divisor, and the second arrow is the adjunction map $\omega_{X/T}(R) \to  i_* i^* \omega_{X/T}(R)= i_* \Omega_{R/T}^1$.

Let $d_{R/T}: \Omega_{R/T}^1 \to \Omega_{R/T}^2$ denote the relative exterior derivative, and set 
\[
\omega_{X/T}' := \ker(\omega_{X/T}(R)  \to i_* \Omega_{R/T}^1 \overset{i_*d_{R/T}}{\longrightarrow} i_*\Omega_{R/T}^2).
\] 
Since both maps are 
$\pinv(\ot)$-linear, $\omega_{X/T}'$ is sheaf of $\pinv(\ot)$-modules on $X$. It is \textbf{not} as $\ox$-module since $i_* d_{R/T}$ is not $\ox$-linear.  We call $\omega_{X/T}'$ the \emph{relative Berezinian sheaf}.

\begin{remark}
    \normalfont For local superconformal coordinates $(z, \theta)$ on $X$ in which $R$ is defined by $z=0$, a global section of $\omega_X'$ is of the form
   \begin{equation} \label{locformomprime} \left (g(z) + zf(z) \theta \right) [dz | d \theta]/z. \end{equation}
\end{remark}

From the definition of $\omega'$, we have a natural map of $\mathcal{O}_T$-modules
\[
\pi_* \omega_{X/T}' \longrightarrow \pi_* i_*\mathcal{Z}_{R/T}^1:=\ker(\pist i_* d_{R/T}).
\]
We denote this morphism by $P_{X/T}^-$ and call it the \emph{odd period map}.




There is also an \emph{even period map}
\[
P_{X/T}^+: \pi_* \omega_{X/T}' \longrightarrow R^1 \pi_* \mathbb{Z} \otimes \mathcal{O}_T,
\qquad
s \longmapsto \left( \gamma \mapsto \int_{\gamma^{th}} s\vert_{X\setminus R}\right)
\]
where we may assume that the $1$-cycle $\gamma \in H_1(X,\z)$ does not meet the Ramond divisor. This is possible because topologically $|R|$ is a finite collection of points.




We can now define the \emph{super period map} to be the direct sum 
\begin{equation}\label{periods}
    P_{X/T}:=P_{X/T}^+\oplus P_{X/T}^-:
    \pi_*\omega_{X/T}'\to\mc{O}_{\Lambda}.
\end{equation} 
where
\[
\mc{O}_{\Lambda}:=
\left(R^1\pist\z\otimes_{\z}\mc{O}_T\right)
\oplus
\pist i_*\mc{Z}_{R/T}^1 
\]

\begin{lemma} \label{super period map is grading-preserving}
The super period map is grading-preserving. 
    
\end{lemma}

\begin{proof}
    It is enough to prove this for $T=\Spec \bb$, where
    \[
    P_X: H^0(X, \omega_X') \to H^1(X, \bb) \oplus H^0(R, \mcz_R^1) \cong \bb^{2g|2r} \qquad s \longmapsto (P_X^+(s), P_X^-(s)).
    \]
    This reduces to proving that $P_X^+(s)=0$ for all $s \in H^0(X, \omega_X')^-$ and $P_X^-(s) = 0$ for all $s \in H^0(X, \omega_X')^+$. 

    Start by choosing local superconformal coordinates $(z, \theta)$, and let $e:=([dz|d \theta]/z)$ denote the local generator of $\omega_X'$. 
    From \eqref{locformomprime}, an even global section of $\omega_X'$ is locally of the form $z \theta f(z) \cdot e$ in these coordinates, and an odd global section is locally of the form $g(z) \cdot e$. 
    It is already immediate that the images of odd global sections under the even period map are all zero. Indeed, this follows from the definition of Berezin integration and the fact that the odd global sections are missing the top $\theta$ term. 

    The odd period map is just restriction to the Ramond divisor. Since $R$ is locally cut out by $z=0$, and since the even sections of $\omega_X'$ are divisible by $z$, they clearly vanish when restricted to $R$.
\end{proof}

For simplicity, set 
\[
\Lambda_0:= R^1 \pist \z, \qquad \Lambda_1= V_{\z}, \qquad \Lambda = \Lambda_0 \oplus \Lambda_1.
\]
For every superscheme $T$, set $\Lambda_T:= \Lambda \otimes_{\z} \ot$. When $T=\Spec \bb$, we write $\Lambda_{\bb}$. 

Let $J_0$ denote the usual symplectic intersection pairing on $\Lambda_0$, recall the quadratic form $J_1$ on $\pist i_* \mcz_{R/T}^1 \cong V_{\z} \otimes_{\z} \ot$ defined in the previous section, and set $J=J_0 \oplus J_1$. Then $J_T:= J \otimes_{\z} \ot$ is a supersymplectic pairing on $\Lambda_{T}$.  

Note that since $\Lambda_{\bb}$ is the complexification of the integral lattice $\Lambda$, complex conjugation is defined:
\[
\overset{-}{{}}: \Lambda_{\bb} \to \Lambda_{\bb}, \qquad \lambda \otimes v \longmapsto \lambda \otimes \ov{v}.
\]
Thus every subspace $V \subset \Lambda_{\bb}$ has a complex conjugate $\ov{V} \subset \Lambda_{\bb}$. 

We say that a subspace $V \subset \Lambda_{\bb}$ is $J$-Lagrangian if $\on{dim}(V)=g|r$ and $J \vert_V \equiv 0$.
For any $V \subset \Lambda_{\bb}$, we always have the induced pairing $J \vert_{V \times \ov{V}}: V \times \ov{V} \to \bb$. If $V$ is $J$-Lagrangian and $V \cap \ov{V}=0$, then $J \vert_{V \times \ov{V}} \to \bb$ is non-degenerate.





\begin{theorem} \label{bigone}
Suppose $X$ is a super Riemann surface with Ramond punctures over $\Spec \bb$ such that
\[
\dim_{\bb} H^0(X,\omega_X')=g|r.
\]
Then the super period map
\[
P_X: H^0(X, \omega_X') \to \Lambda_{\bb}:=H^1(X, \bb) \oplus H^0(R, \mcz_R^1), \qquad s \longmapsto (P_X^+(s), P_X^-(s))
\]
is injective, and its image $\on{im}(P_X)$ is a $J$-Lagrangian subspace of $\Lambda_{\bb}$. 
Furthermore, $\on{im}(P_X) \cap \ov{\on{im}(P_X)}=0$, and thus $J \vert_{\on{im}(P_X) \times \ov{\on{im}(P_X)}}$ is non-degenerate. 

\end{theorem}

\begin{proof}
The main ingredient of the proof is the decomposition 
\begin{equation} \label{decomp of ext ber in split case}
H^0(X, \omega_X') = H^0(C, \omega_C) \oplus \Pi H^0(C,L).
\end{equation}

The kernel of $P_X$ is the intersection $\ker(P_X^+) \cap \ker(P_X^-)$. Since $P_X$ is grading-preserving,
\[
\Pi H^0(C, L) \subseteq \ker(P_X^+), \qquad H^0(C, \omega_C) \subseteq \ker(P_X^-).
\]
Both of these inclusions are equalities. Indeed, the restriction of $P_X^+$ to $H^0(C, \omega_C)$ is just the ordinary period map, which is well-known to be injective. 
To prove the claim for $P_X^-$, we use the identification $H^0(R, \mcz_R^1)=\Pi L$ from Lemma \ref{oddtarget}, and note that the restriction of $P_X^-$ to $\Pi H^0(C,L)$ is just the map $H^0(C, L) \to L \vert_D$ in the sequence
\begin{equation} \label{restriction sequence}
0 \longrightarrow H^0(C,L(-D)) \longrightarrow H^0(C,L) \overset{\fermp}{\longrightarrow} L \vert_D \to H^0(C,L)^{\vee} \longrightarrow H^0(C,L(-D))^{\vee} \longrightarrow 0.
\end{equation}

Since $H^0(X,\omega_X')=g|r$, $H^1(C,L)=H^0(C,L(-D))=0$, and hence $\fermp$ is injective.

For the second claim, the image of $P_X$ in $\Lambda_{\bb}$ is said to be $J$-Lagrangian if the respective images of $P_X^+$ and $P_X^-$ in $H^1(C,\bb)$ and $H^0(R, \mcz_R^1)$ are $J_0$- and $J_1$-maximally isotropic, respectively. 
The image of $P_X+$ is $J_0$-maximally isotropic since it coincides with the image of the ordinary period map for $C$, which is well-known to be maximally isotropic with respect to the symplectic intersection pairing $J_0$.

The claim for the image of the odd period map follows from the identification of its image with the image of $H^0(C,L) \to \Pi L \vert_D$, identification of $J_1$ with the residue pairing on $L \vert_D$ (c.f. Lemma \ref{oddtarget}), and the residue theorem:  Since $s \otimes s \in H^0(C,L^2)=H^0(C, \omega_C(D))$ is a global meromorphic one-form on the curve $C$, $\sum_{i=1}^{2r} \on{res}_{p_i}(s \otimes s)=0$.  

For the remaining claim, we first note that the target $\Lambda_{\bb}:=\Lambda \otimes \bb$ comes with a complex-conjugation map: $\lambda \otimes v \mapsto \lambda \otimes \ov{v}$. Thus, we can define the complex-conjugate super period map $\ov{P}_X(\ov{s}):= \ov{P_x(s)}$ for all $s \in H^0(X, \omega_X')$, and note that
\[
\ov{\on{im}(P_X)} = \on{im}(\ov{P}_X).
\]
To prove the transversality statement, it suffices to show that
\[
\on{im}(P_X^+) \cap \on{im}(\ov{P_X^+})=0, \qquad \on{im}(P_X^-) \cap \on{im}(\ov{P_X^-}) =0.
\]
The claim for the even period map again follows from its identification with the image of the ordinary period map, which is well-known to be transverse to its complex conjugate. To prove the claim for the odd period map, assume the intersection is non-zero. Then there exists $v \in \Lambda_{\bb}$ such that $v= \ov{v}$. 
But then $v$ is a nonzero real isotropic vector for $J_1$, which is impossible since $J_1$ is positive definite. 
\end{proof}




\section{The bad locus} \label{sectiontitle: baddys}

Throughout this section we work entirely on the bosonic reduction.  Let
\(M_{\bos}\) denote the bosonic reduction of the chosen \'etale chart
\(M\to\mfr\), and let \(B_{\bos}\subset M_{\bos}\) denote the bosonic bad
locus.  A point of \(M_{\bos}\) is the twisted-spin datum
\[
(C,D,L),\qquad D=p_1+\cdots+p_{2r},\qquad L^2\simeq K_C(D).
\]
Although Ramond punctures cannot be forgotten at the level of super Riemann
surfaces, there is an ordinary bosonic morphism
\[
p:M_{\bos}\longrightarrow \mc M_g,\qquad (C,D,L)\longmapsto C.
\]
Theorems \ref{dimensioncounts} and \ref{exactcodim} below are statements about
this bosonic morphism and the ordinary complex codimension of
\(B_{\bos}\subset M_{\bos}\).  Equivalently, they give the even part of the
super-codimension of the corresponding closed locus in \(M\).  For readability
we suppress the subscript \(\bos\) inside the proofs.

\begin{theorem} \label{dimensioncounts}
If $r>1$, then every irreducible component of the bosonic bad locus
$B_{\bos}$ has complex codimension $>1$ in $M_{\bos}$.
\end{theorem}

\begin{proof}
On the bosonic reduction the bad locus is
\[
B=\{(C,D,L):h^0(C,L)\ge r+1\}
\subset
M=\{(C,D,L):L^2\simeq K_C(D)\}.
\]
Set
\[
A:=K_C\otimes L^{-1}=L(-D).
\]
Then
\[
A^2(D)=K_C,\qquad \deg A=g-1-r,
\]
and Riemann--Roch gives
\[
h^0(A)=h^0(L)-r.
\]
Thus the bad condition is $h^0(A)>0$.

Fix a curve $C$.  Write
\[
F:\operatorname{Pic}^{g-1-r}(C)\longrightarrow \operatorname{Pic}^{2r}(C),
\qquad F(A):=K_C\otimes A^{-2}.
\]
Then the fiber of $M\to\mc M_g$ over $C$ is
\[
M_C=
F^{-1}(W^0_{2r})\times_{W^0_{2r}}\operatorname{Sym}^{2r}(C).
\]
Equivalently, the map
\[
\rho_C:M_C\longrightarrow \operatorname{Sym}^{2r}(C),\qquad (D,L)\longmapsto D,
\]
is finite: for every effective divisor $D$ there are exactly $2^{2g}$ square
roots of $K_C(D)$.  In particular,
\[
\dim M_C=2r
\]
for every $C$.

Assume, for a contradiction, that $B$ has an irreducible component $B'$ of
codimension one.  Let
\[
Y=\overline{p(B')}\subset\mc M_g,
\qquad c=\operatorname{codim}_{\mc M_g}Y.
\]
Since $\dim M=3g-3+2r$, the generic fiber of $B'\to Y$ has dimension
\[
(3g-4+2r)-(3g-3-c)=2r+c-1.
\]
This cannot exceed $\dim M_C=2r$, hence $c\le1$.  There are therefore two
possibilities.

\smallskip
\noindent
\emph{Case (1): $c=0$.}
Then $B'$ dominates $\mc M_g$, and for a general curve $C$ the fiber
\[
B'_C:=B'\cap M_C
\]
is a divisor in $M_C$.  Since $\rho_C$ is finite, its image
\[
Z_C:=\rho_C(B'_C)\subset\operatorname{Sym}^{2r}(C)
\]
is a divisor.  Consider the closed $r$-dimensional subvariety
\[
Y_C:=\{2E:E\in\operatorname{Sym}^r(C)\}\subset\operatorname{Sym}^{2r}(C).
\]
The symmetric product $\operatorname{Sym}^{2r}(C)$ is projective of dimension $2r$, so the
projective dimension theorem gives
\[
\dim(Z_C\cap Y_C)\ge(2r-1)+r-2r=r-1.
\]
For every $2E\in Z_C\cap Y_C$ choose a bad lift $(C,2E,L)\in B'_C$ and set
\[
N:=L(-E).
\]
Then $N^2\simeq K_C$, while
\[
L(-2E)=N(-E)
\]
has a nonzero section.  Hence, for an $(r-1)$-dimensional family of effective
divisors $E$ of degree $r$, the divisor $E$ is contained in the zero divisor
of a section of some theta characteristic $N$.  Since a fixed curve has only
finitely many theta characteristics, one fixed $N$ works for an
$(r-1)$-dimensional family of $E$'s.  For fixed $N$, the space of such $E$'s
has dimension at most $\dim|N|=h^0(N)-1$, because a fixed divisor in $|N|$
has only finitely many degree-$r$ effective subdivisors.  Thus
\[
h^0(N)\ge r.
\]
For $r\ge2$ this contradicts the fact that on a general curve every theta
characteristic has at most one section.

\smallskip
\noindent
\emph{Case (2): $c=1$.}
Then $Y\subset\mc M_g$ is a divisor, and for a general $C\in Y$ the fiber
$B'_C$ contains an irreducible component of $M_C$.  Every irreducible
component of $M_C$ maps surjectively to $\operatorname{Sym}^{2r}(C)$ under the finite map
$\rho_C$.  Hence for every $E\in\operatorname{Sym}^r(C)$ there is a bad lift with
$D=2E$.  As above, writing $N=L(-E)$ produces a theta characteristic for
which $N(-E)$ is effective.  Since there are finitely many theta
characteristics, one fixed $N$ works for an $r$-dimensional family of $E$'s,
and therefore
\[
h^0(N)\ge r+1.
\]
For $r\ge2$, Corollary 1.3 of \cite{polishchuk2003effective}, specialized to
spin order $2$ with no marked points, shows that the locus of curves carrying
a theta characteristic $N$ with $h^0(N)\ge r+1$ has codimension at least
$r$ in $\mc M_g$.  It therefore cannot contain the divisor $Y$, a
contradiction.

Both cases are impossible, proving the theorem.
\end{proof}

The following theorem strengthens Theorem \ref{dimensioncounts}; its proof uses Polishchuk's moduli space of effective spin structures.

\begin{theorem} \label{exactcodim}
Assume $r>0$. If $1\leq r\leq g-1$, every irreducible component of the
bosonic bad locus $B_{\bos}\subset M_{\bos}$ has codimension exactly $r$.
If $r>g-1$, the bosonic bad locus is empty.
\end{theorem}

\begin{proof}
This is again a purely bosonic statement.  If $r>g-1$, then
\[
\deg A=g-1-r<0,
\]
so $h^0(C,A)=0$ and the bad locus is empty.  We therefore assume
\[
1\leq r\leq g-1.
\]

It is convenient first to label the $2r$ points of the Ramond divisor.
Consider Polishchuk's ordinary moduli stack of effective spin structures
\cite{polishchuk2003effective} with spin order $2$ and parameters
\[
n=2r,\qquad m_1=\cdots=m_{2r}=1,\qquad d=g-1-r.
\]
Thus this stack, which we denote by $\mc E$, parametrizes data
\[
(C,E,p_1,\ldots,p_{2r}),\qquad \deg E=d,
\]
satisfying
\begin{equation} \label{effective spin equation}
\mc O_C(2E+p_1+\cdots+p_{2r})\simeq K_C.
\end{equation}
Polishchuk's Theorem 1.1(a) says that $\mc E$ is smooth of dimension
\[
2g-2+d+n
=2g-2+(g-1-r)+2r
=3g-3+r.
\]

We also need that $\mc E$ is irreducible. By Corollaries 4.1 and 4.2 of
\cite{polishchuk2003effective}, it is connected in our case. Indeed, for
$g\geq4$ none of the exceptional disconnected cases occurs, since all the
marked multiplicities $m_i$ are equal to $1$; for $g=2$ all the relevant
spaces are connected; and for $g=3$ none of the exceptional cases has
$n=2r$ and all $m_i=1$. Since $\mc E$ is smooth, its irreducible components
are disjoint open-and-closed substacks. Hence connectedness implies that
$\mc E$ is irreducible.

Let $\mc S$ be the ordinary moduli stack of the corresponding labeled
$2$-spin data
\[
(C,p_1,\ldots,p_{2r},A),\qquad
A^{\otimes2}\simeq K_C(-p_1-\cdots-p_{2r}),
\]
and let $B_{\rm lab}\subset\mc S$ be the locus where $h^0(C,A)>0$.
Polishchuk points out that the map from effective spin structures to the
spin stack is only locally defined on the stack because of scalar
automorphisms of the spin structure, but it is an honest morphism on coarse
moduli spaces.  We therefore make the next dimension argument on coarse
moduli spaces.  For these Deligne--Mumford stacks over $\bb$, passage to the
coarse moduli space preserves dimensions and irreducible components.

On coarse moduli spaces we have a surjective morphism
\[
\Phi:\mc E_{\rm coarse}\longrightarrow (B_{\rm lab})_{\rm coarse},\qquad
(C,E,p_1,\ldots,p_{2r})\longmapsto
(C,p_1,\ldots,p_{2r},A=\mc O_C(E)).
\]
Equation \eqref{effective spin equation} is precisely the condition
$A^{\otimes2}\simeq K_C(-D)$, where $D=p_1+\cdots+p_{2r}$. Conversely, if
$A$ has a nonzero section, the zero divisor of that section gives such an
$E$, proving surjectivity.

The geometric fiber over $(C,p_1,\ldots,p_{2r},A)$ is a finite quotient of
the complete linear system $\mathbb P H^0(C,A)$ by the finite automorphism
group of the datum; in particular its dimension is $h^0(C,A)-1$. By Theorem
1.1(b) of \cite{polishchuk2003effective}, there exists a point of $\mc E$ for
which
\[
h^0(C,\mc O_C(E))=1.
\]
Upper semicontinuity therefore shows that the locus where
$h^0(C,\mc O_C(E))=1$ is a nonempty open subset of the irreducible stack
$\mc E$, hence is dense. Thus $\Phi$ is generically finite. Consequently
\[
\dim B_{\rm lab}=\dim\mc E=3g-3+r.
\]
Since the coarse space of $\mc E$ is irreducible, its image
$(B_{\rm lab})_{\rm coarse}$ is irreducible, and hence so is
$B_{\rm lab}$.

The ordinary stack $\mc S$ of labeled $2$-spin data is finite over
$\mc M_{g,2r}$ and hence has dimension
\[
3g-3+2r.
\]
It follows that
\[
\operatorname{codim}_{\mc S} B_{\rm lab}
=(3g-3+2r)-(3g-3+r)=r.
\]
Finally, forgetting the labels on the Ramond points is a finite quotient, and
pulling back to the bosonic reduction $M_{\bos}$ of the \'etale chart used in
this paper preserves codimension. Therefore every irreducible component of
$B_{\bos}\subset M_{\bos}$ has codimension exactly $r$.
\end{proof}

The next proposition is not needed for the main theorem, but we include it for completeness.

\begin{proposition}
    The bad locus $\mc{B} \subset M$ descends to a closed substack of $\mathfrak{M}_{g,0,2r}$ .
\end{proposition}

\begin{proof}
Let $S \to \mc{S}\mc{M}_{g,0,2r}$ denote the restriction of $M \to \mfr$ to the twisted spin moduli space $\mc{S}\mc{M}_{g,0,2r} \subset \mfr$, and note that since $|M|=|S|$, the bad locus $\mc{B}$ is also naturally a subset of $|S|$.
Since $\mc{S} \mc{M}_{g,0,2r} \subset \mathfrak{M}_{g,2r}$ is a closed substack, it is enough to prove that $\mc{B} \subset |S|$ descends to a closed substack of $\mc{S} \mc{M}_{g,0,2r}$  \footnote{The statement is well-known for $r=0$ where the bad locus in the spin moduli space coincides the the theta-null divisor.}.

Let $\pi: C \to S$ denote the pullback of the universal spin curve, with Ramond divisor $R$ and twisted spin structure $L$.
Choose  a relative effective Cartier divisor $D \subset C$ of degree $d \gg 0$ such that $H^1(C_s, L_s(D_s))=0$ for all $s \in S$. Pushing forward $0 \to L \to L(D) \to L(D)\vert_D \to 0$ to $S$ gives
 \[ 0 \to \pist L \to \pist L(D) \overset{\alpha}{\to} \pist L(D) \vert_D \to R^1 \pist L \to 0 \]
where $\pist L(D)$ and $\pist L(D) \vert_D$ are locally free of rank $r+d$ and $d$, respectively. After restriction to $s \in S$, one has $\ker(\alpha_s) = H^0(C_s, L_s)$.  Hence $H^0(C_s, L_s) >r$ when $\alpha_s$ is not surjective.

The bad locus on $S$ is therefore the closed subscheme  $V(\on{Fitt}_0(\on{coker}(\alpha)))$, where $\on{Fitt}$ denotes the Fitting ideal; this is the invariant ideal which is locally generated by the maximal minors of a matrix representing $\alpha$.
Finally, this closed subscheme descends to a closed substack of $\mc{S} \mc{M}_{g,2r}$ because its defining ideal commutes with base change, and so the pullbacks of $\mc{B}$ to $S \times_{\mc{S} \mc{M}_{g,2r}} S$ agree.
\end{proof}

\section{Berezinian vs. extended Berezinian}

Continue to let $M \to \mfr$ denote an \'etale cover, and $\pi : X \to M$ the pullback of the universal family.

For simplicity set $\omega=\omega_{X/M}$, $\omega'=\omega_{X/M}'$, and $Q=i_*\mc Z_{R/M}^1$. Since $\omega$ is the kernel of the canonical map $\omega' \to i_* \mc{Z}_{R/M}^1$ we have a short exact sequence 
\[ 0 \to \omega \xrightarrow{\iota} \omega' \to Q \to 0\]
of $\pinv(\mc{O}_M)$-modules. This gives an exact triangle
\[ R \pist \omega \to R \pist \omega' \to R \pist Q\]
in the derived category $D(\mc O_M)$, and an isomorphism 
\[ \Ber(R \pist \omega') \cong   \Ber(R \pist \omega ) \otimes \Ber(R\pist Q). \]

Fix an orientation $\sigma$ on $R$. This determines a trivialization 
\[ \pist i_* \mcz_{R/M}^1 \cong \mc{O}_M^{0|2r}. \]
Therefore: 

\begin{theorem}\label{firstresult} Suppose $r>0$. After fixing an orientation $\sigma$ on the Ramond divisor $R$, there is a canonical isomorphism
\[ \Ber(R \pist \omega_{\mcx/M}) \cong \Ber(R \pist \omega_{\mcx/M}')\]
    
\end{theorem}

\begin{lemma}\label{orientation character}
Let $\sigma$ and $\sigma'$ be two orientations related by
\[
 dx'^i=\epsilon_i dx^i,\qquad \epsilon_i\in\{\pm1\},
 \qquad
 \chi(\epsilon):=\prod_{i=1}^{2r}\epsilon_i.
\]
Let $\alpha_{\sigma}$ and $\alpha_{\sigma'}$ denote the two isomorphisms of Theorem \ref{firstresult}. Then
\[
 \alpha_{\sigma'}=\chi(\epsilon)\,\alpha_{\sigma}.
\]
\end{lemma}

\begin{proof}
The quotient $Q=i_*\mc Z^1_{R/M}$ has no higher direct images, so an orientation identifies
\[
R\pi_*Q=\pi_*Q\cong\mc O_M^{0|2r}.
\]
Changing the orientation by $\epsilon$ changes this purely odd frame by the diagonal matrix
\[
S_{\epsilon}=\operatorname{diag}(\epsilon_1,\ldots,\epsilon_{2r}).
\]
The induced change of the Berezinian generator is
\[
\Ber(S_{\epsilon})=\det(S_{\epsilon})^{-1}
=\chi(\epsilon),
\]
since $\chi(\epsilon)=\pm1$.  The isomorphism of Theorem \ref{firstresult} is obtained from the exact triangle
\[
R\pi_*\omega\longrightarrow R\pi_*\omega'
\longrightarrow R\pi_*Q
\]
by trivializing $\Ber(R\pi_*Q)$ with this generator.  Hence it changes by the same factor $\chi(\epsilon)$.
\end{proof}

 \section{Measure over the good locus} \label{analysis of super period map on good locus}\label{sectiontitle: measure over good locus}

Throughout this section, we fix an \'etale cover $M \to \mfr$ and let $\pi: \mcx \to M$ denote the pullback of the universal family to $M$. We work on the analytification of this smooth chart, and hence regard $M$ as a complex supermanifold. No properness of the smooth moduli space is used.

Let $\ov{M}=(|M|, \mc{O}_{\ovm})$ denote the complex-conjugate supermanifold of $M$, where $\mc{O}_{\ov{M}}$ is the sheaf of anti-holomorphic functions on $M$. Let $M_{\rr}$ denote the smooth supermanifold underlying $M$. There is a natural map $\Delta:=(\on{id}, c): M_{\rr} \to M \times \ov{M}$ sending $x \mapsto (x,x)$ and satisfying
\[
\deltpull F(z(x), \theta(x), w(x), \eta(x)) = F(z(x), \theta(x), \ov{z(x)}, \ov{\eta(x)}).
\]
Let $p_1: M \times \ov{M} \to M$ and $p_2: M \times \ov{M} \to \ov{M}$ denote the canonical projections. For every subset $Z \subset \mrr$, let $\Delta_Z: Z \to \mrr \xra{\Delta} M \times \ovm$ denote the inclusion. For every coherent sheaf $E$ on $M$ and every $Z \subset \mrr$, there are natural smooth isomorphisms
\[
\deltpull_Z \pons E \cong E \vert_Z, \qquad \deltpull_Z \pows \ov{E} \cong \ov{E} \vert_Z.
\]
If $Z=\{z\}$ is a point, then these are isomorphisms of complex super vector spaces.

The complex conjugate of a complex supermanifold is discussed in Section 5.1 of \cite{felder2019regularity} and Section 4.6 of \cite{deligne1999notes}.

For simplicity, set $P:=P_{\mcx/M}$ to denote the super period map. Choose an orientation of the Ramond divisor temporarily. Since the local systems involved are locally constant, after shrinking to a sufficiently small open neighborhood $\mc N$ of the diagonal in $M\times\ov M$ we may identify the two pullbacks of $\Lambda$ compatibly with the pairing $J$. We henceforth work on this neighborhood. Let $\beta$ denote the holomorphic pairing
\begin{equation} \label{bilinear pairing}
\pons \pist \omprime \otimes \pows \pist \ov{\omprime}
\xra{\pons P_{\mcx/M} \otimes \pows \ov{P_{\mcx/M}}}
(\Lambda \otimes_{\z} \Lambda) \otimes_{\z} \mc{O}_{\mc N}
\xra{J}
\mc{O}_{\mc N},
\end{equation}
sending $(f, \ov{g}) \longmapsto J(P(f), \ov{P}(\ov{g}))=J(P(f), \ov{P(g)})$. For every $Z \subset \mrr$, set $\beta_Z:=\deltpull_Z \beta$ whenever the corresponding part of the diagonal lies in $\mc N$.

\begin{lemma}\label{good locus base change}
On the good locus $U$, the sheaf $\pi_*\omega'$ is locally free of rank
$g|r$, the sheaf $R^1\pi_*\omega'$ is locally free of rank $1|0$, and
formation of both sheaves commutes with restriction to fibers (more generally,
with base change inside $U$). Moreover, $R^1\pi_*\omega'$ has a canonical
trivialization $R^1\pi_*\omega'\cong\mc O_U$.
\end{lemma}

\begin{proof}
The exact sequence
\[
0\longrightarrow\omega\longrightarrow\omega'\longrightarrow Q\longrightarrow0
\]
gives an exact triangle after applying $R\pi_*$.  The morphism $\pi$ is proper
and smooth and $\omega$ is locally free, so $R\pi_*\omega$ is perfect.  The
sheaf $Q$ is supported on the finite Ramond divisor, and $R\pi_*Q=\pi_*Q$ is
locally free and concentrated in degree zero.  Hence $R\pi_*\omega'$ is also
perfect, of amplitude contained in $[0,1]$.

On $U$ the fiber dimensions in \eqref{dims} are constant:
\[
h^0(\omega'_X)=g|r,\qquad h^1(\omega'_X)=1|0.
\]
Locally represent $R\pi_*\omega'$ by a two-term complex of finite free
modules.  Constancy of these two cohomology ranks is equivalent to constancy
of the rank of its differential.  Its kernel and cokernel are therefore
locally free and commute with base change.  These are respectively
$\pi_*\omega'$ and $R^1\pi_*\omega'$.  In particular, restriction to a closed
point $x\in U$ identifies $\pi_*\omega'|_x$ with $H^0(X_x,\omega'_{X_x})$.

For completeness, the trivialization of $R^1\pi_*\omega'$ is also canonical.
The Serre-duality statement for the extended Berezinian identifies
$R^1\pi_*\omega'$ with the dual of the sheaf of holomorphic functions whose
restriction to each Ramond divisor is constant; see \cite[Section~5.1]{witten2015super}.
On the good locus this latter sheaf is canonically $\mc O_U$, generated by the
constant function $1$.  Dualizing gives
\[
R^1\pi_*\omega'\cong\mc O_U
\]
canonically.
\end{proof}

\begin{proposition} \label{lagrangian image of the good locus} The restriction of $P$ to the good locus $U \subset M$ is injective, and its image is a $J_U$-Lagrangian subbundle of $\Lambda_U$. 
    
\end{proposition}

\begin{proof}
By Lemma \ref{good locus base change}, $\pi_*\omega'|_U$ is a vector bundle and its restriction to every closed point $x\in U$ is $H^0(X_x,\omega'_{X_x})$. The period map is compatible with this restriction. Hence injectivity and the Lagrangian condition may be checked on closed points, where they are exactly Theorem \ref{bigone}.
\end{proof}

\begin{corollary} \label{corollary of big one}
Fix an orientation $\sigma$ on $R$. The restriction of the pairing \eqref{bilinear pairing} to the good locus $U_{\rr} \subset \mrr$ is non-degenerate. In particular, combining this with Theorem \ref{firstresult} and the super Mumford isomorphism, we get a nowhere-vanishing smooth section
\[
s_U \in \Gamma(U_{\rr}, \deltpull_U(\pons \Ber^5(R\pist \omega_{\mcx/M}) \otimes \pows \Ber^5(R\ov{\pist \omega_{\mcx/M}}))).
\]
\end{corollary}

\begin{proof}
Over the good locus in $\mrr$, the target of the pairing
\[
\beta_U: \pist\omega' \vertu \otimes \ov{\pist \omega'} \vertu
\xra{\deltpull_U \pons P \otimes \deltpull_U \pows\ov{P}}
\Lambda_{\urr} \otimes \Lambda_{\urr}
\xra{J_{\urr}}
\mc{O}_{U_{\rr}}
\]
is a smooth vector bundle, so it suffices to check non-degeneracy after restriction to closed points. For every  closed point $x \in U$, we have canonical isomorphisms of super vector spaces:
\[
\deltpull_x \pons \omega' \cong \pist \omega' \vert_x, \qquad
\deltpull_x \pows \ov{\pist \omega'} \cong \ov{\pist \omega'} \vert_x, \qquad
\deltpull_x \mc{O}_{M \times \ovm} \cong \bb.
\]
Since $\pist \omega'$ and its complex conjugate are locally free in a neighborhood of $x$, the base change maps give natural isomorphisms
\[
\pist \omega' \vert_x \cong H^0(X, \omega_X'), \qquad
\ov{\pist \omega'} \vert_x \cong \ov{H^0(X, \omega_X')},
\]
where $X$ denotes the fiber $\mcx_x$. In particular, $\deltpull_x \pons P = P_X$ and $\deltpull_x \pows \ov{P} = \ov{P_X}$, and thus
\[
\beta_x: H^0(X, \omega_X') \otimes \ov{H^0(X, \omega_X')}
\to \Lambda_{\bb} \otimes \Lambda_{\bb}
\xra{J_{\bb}}
\bb,
\]
which is non-degenerate by Theorem \ref{bigone}.

For the second claim, take $\Ber$ and apply Theorem \ref{firstresult} and the super Mumford isomorphism. 
\end{proof}

\begin{corollary}\label{orientation independent measure}
The measure $s_U$ of Corollary \ref{corollary of big one} is independent of the choice of orientation of the Ramond divisor.  In particular, it descends from the oriented cover to the original moduli space.
\end{corollary}

\begin{proof}
Let $\sigma'$ be obtained from $\sigma$ by
$dx^i\mapsto \epsilon_i dx^i$, and put
$S_{\epsilon}=\operatorname{diag}(\epsilon_1,\ldots,\epsilon_{2r})$ and
$\chi(\epsilon)=\prod_i\epsilon_i$.  The odd period vectors are transformed by $S_{\epsilon}$, while
\[
S_{\epsilon}^{t}J_1S_{\epsilon}=J_1.
\]
Thus the period pairing $\beta$ is unchanged when the orientation is changed simultaneously on the holomorphic factor and on its complex conjugate.

By Lemma \ref{orientation character}, the isomorphism of Theorem \ref{firstresult} changes by $\chi(\epsilon)$.  After taking the fifth power, the holomorphic factor therefore changes by
$\chi(\epsilon)^5=\chi(\epsilon)$, and the conjugate factor changes by the same sign.  Their tensor product is multiplied by
\[
\chi(\epsilon)^2=1.
\]
Hence the resulting section $s_U$ is unchanged.

The same calculation applies holomorphically on the neighborhood $\mc N$ of
the diagonal: identify the two orientation torsors using the chosen
identification of the two copies of $\Lambda_1$, choose an orientation on one
factor, and transport it to the other.  A simultaneous change by $\epsilon$
leaves the period pairing invariant and multiplies the two fifth-power
Berezinian identifications by the same character $\chi(\epsilon)$, so their
product is unchanged.  Thus the determinant-level construction used below is
independent of this paired local orientation already before restricting to the
smooth diagonal.
\end{proof}




\begin{remark}
We do not see a reasonable extension of this result to a neighborhood of the quasidiagonal. The problem arises from Theorem \ref{bigone}: on the diagonal we compare the image of $P_{-}$ to its conjugate and we know that these are transversal. But on the quasidiagonal, where we allow $L'$ to be independent of $L$, we lose all control over $P_{-}$.
\end{remark}
\section{Extension of the measure across the bad locus} \label{sectiontitle: extension across the bad locus}

Let $B_{\bos}\subset M_{\bos}$ be the bosonic bad locus of Section
\ref{sectiontitle: baddys}, and let $U\subset M$ be the open subsupermanifold
whose bosonic reduction is $M_{\bos}\setminus B_{\bos}$.  Thus every
codimension statement used below is a statement on the bosonic reduction.
We work in a sufficiently small open neighborhood $\mc N$ of the diagonal in $M\times\ov M$ on which the two copies of the local system $\Lambda$ are identified, as in Section \ref{sectiontitle: measure over good locus}. This identification also identifies the two orientation torsors. We may therefore choose a local orientation on one factor and transport it to the other. By the final paragraph of Corollary \ref{orientation independent measure}, all determinant-level objects below are independent of this paired auxiliary choice and hence descend to the unoriented neighborhood. Let
\[
\mc{L}:=p_1^*\Ber^5(R\pi_*\omega)\otimes p_2^*\Ber^5(R\ov{\pi}_*\ov{\omega}).
\]
This is a holomorphic line bundle on $\mc{N}$. Over
$(U\times \ov U)\cap \mc{N}$, Theorem \ref{firstresult} and the canonical
trivialization of $R^1\pi_*\omega'$ identify $\mc L$ with
\[
\Ber^5(E_1)\otimes \Ber^5(E_2),
\qquad
E_1:=p_1^*\pi_*\omega',\quad
E_2:=p_2^*\ov{\pi_*\omega'}.
\]
Let $\beta:E_1\otimes E_2\to\mc O$ be the period pairing of Section
\ref{sectiontitle: measure over good locus}. Its non-degeneracy locus contains
$\Delta(U_{\rr})$. On that locus the inverse fifth Berezinian of $\beta$ is
the holomorphic section whose pullback to $U_{\rr}$ is the measure constructed
in Corollary \ref{corollary of big one}.

\begin{proposition}\label{direct extension proposition}
Assume $r>1$. Let $b\in B_{\bos}$ and let $V\ni b$ be a sufficiently small
open neighborhood in $M$. After shrinking $V$ if necessary, the measure
section on the good locus extends holomorphically to a section
\[
\widehat t\in \Gamma\bigl(V\times \ov V,\mc L\bigr).
\]
Consequently, its pullback by the smooth diagonal extends the smooth measure
across $b$. We do not claim that this extension is nowhere vanishing.
\end{proposition}

\begin{proof}
Shrink $V$ so that
\[
V\times\ov V\subset\mc N
\]
and so that the reduced complex manifold
$V_{\bos}\setminus B_{\bos}$ is connected.  Set
\[
W:=(V\cap U)\times(\ov V\cap\ov U),
\]
and let $W^{\times}\subset W$ be the locus where $\beta$ is non-degenerate.
By Corollary \ref{corollary of big one}, $W^{\times}$ contains the portion of
the smooth diagonal lying over $V\cap U$.

After shrinking $V$ further, the even-even part of $\beta$ is non-degenerate
on all of $W$. Its bosonic reduction is the ordinary Hodge pairing, which is
non-degenerate on the conjugate diagonal; hence its determinant remains a unit
on a sufficiently small neighborhood. Cover $W$ by open sets $W_i$ on which
$E_1$ and $E_2$ are trivial and choose local frames. These frames determine a
local frame $e_i$ of $\mc L|_{W_i}$. Write the matrix of $\beta$ as
\[
G_i=
\begin{pmatrix}
A_i&B_i\\
C_i&D_i
\end{pmatrix}.
\]
The matrix $A_i$ is invertible on all of $W_i$.  Therefore the full pairing is
non-degenerate precisely where the holomorphic Schur-complement determinant
\[
\delta_i:=\det(D_i-C_iA_i^{-1}B_i)
\]
is nonzero.  Since $W_{\bos}$ is connected and $W^{\times}$ is nonempty, the
analytic nonvanishing locus $W^{\times}_{\bos}$ is dense in $W_{\bos}$.
Consequently $(W_i\cap W_j)_{\bos}\cap W^{\times}_{\bos}$ is dense in every
nonempty overlap $(W_i\cap W_j)_{\bos}$.

On $W_i\cap W^{\times}$ we have
\[
q_i:=\Ber(G_i)^{-1}
=
\frac{\det(D_i-C_iA_i^{-1}B_i)}{\det A_i}.
\]
The expression on the right is holomorphic on all of $W_i$, whether or not
$G_i$ is invertible. Denote this holomorphic extension by $\widetilde q_i$,
and set
\[
\widetilde t_i:=\widetilde q_i^{\,5}e_i\in\Gamma(W_i,\mc L).
\]

We check that these local sections glue. Let $S_{ij}$ and $T_{ij}$ be the
change-of-frame matrices for $E_1$ and $E_2$, and put
\[
a_{ij}:=\Ber(S_{ij})\Ber(T_{ij})\in\mc O(W_i\cap W_j)^\times.
\]
Then
\[
e_j=a_{ij}^{5}e_i,
\qquad
q_j=a_{ij}^{-1}q_i
\]
on $(W_i\cap W_j)\cap W^{\times}$. Hence the holomorphic superfunction
\[
\widetilde q_j-a_{ij}^{-1}\widetilde q_i
\]
vanishes on a dense open subset of the reduced overlap. Expanding in local odd
coordinates, every coefficient vanishes identically, and therefore
\[
\widetilde q_j=a_{ij}^{-1}\widetilde q_i
\]
on $W_i\cap W_j$. Consequently the $\widetilde t_i$ glue to a holomorphic
section
\[
\widetilde t\in\Gamma(W,\mc L)
\]
which restricts to the original measure section on $W^{\times}$.

It remains to extend $\widetilde t$ across the factors of the bad locus. On
the bosonic reduction,
\[
(V\times\ov V)_{\bos}\setminus W_{\bos}
=
((V_{\bos}\cap B_{\bos})\times\ov V_{\bos})
\cup
(V_{\bos}\times(\ov V_{\bos}\cap\ov B_{\bos})).
\]
By Theorem \ref{dimensioncounts}, each component of this analytic subset has
complex codimension at least $2$. Trivialize $\mc L$ on supercoordinate
neighborhoods $Q_a\subset V\times\ov V$ and write
\[
\widetilde t|_{Q_a\cap W}=f_as_a.
\]
Ordinary Hartogs' theorem on $(Q_a)_{\bos}$, applied coefficientwise after
expanding $f_a$ in the local odd coordinates, extends $f_a$ uniquely to a
holomorphic superfunction $\widetilde f_a$ on $Q_a$. Uniqueness shows that
the sections $\widetilde f_as_a$ agree on overlaps. They therefore glue to the
required section $\widehat t$.
\end{proof}

\begin{theorem}\label{main}
Assume $g\ge2$ and $r>1$. The measure constructed on the good locus extends holomorphically to an open neighborhood of the entire diagonal in $M\times\ov M$. Its pullback to $M_{\rr}$ is therefore a smooth measure extending the measure of Corollary \ref{corollary of big one}. The resulting measure is independent of the auxiliary choice of orientation of the Ramond divisor.
\end{theorem}

\begin{proof}
Apply Proposition \ref{direct extension proposition} near every point of the bad locus. On overlaps the local extensions agree, because they agree on the dense good locus and holomorphic sections of a line bundle which agree on a dense open subset agree everywhere. Thus the local extensions glue to a holomorphic section on a neighborhood of the whole diagonal. Pullback by $\Delta:M_{\rr}\to M\times\ov M$ gives the required smooth measure.

By Corollary \ref{orientation independent measure}, two choices of orientation give the same section on the dense good locus. Their holomorphic extensions therefore agree everywhere by uniqueness. Thus the extended measure is orientation-independent as well.
\end{proof}

\begin{remark}[Descent to supermoduli]
The preceding construction was carried out on an \'etale chart
$M\to\mfr$. It is compatible with \'etale base change and therefore
descends to $\mfr$. Indeed, on
\[
M^{(2)}:=M\times_{\mfr}M
\]
the two pullbacks of the measure agree over the inverse image of the good
locus, since there they are obtained from the same canonical period pairing
and super Mumford isomorphism. This inverse image is dense, and both sections
are holomorphic extensions of their common restriction. By uniqueness of
holomorphic extension they therefore agree on all of $M^{(2)}$. The usual
descent condition follows, so the resulting measure is a well-defined section
on the moduli superstack itself.
\end{remark}

The preceding theorem deliberately asserts extension as a measure, not as a volume form. In fact, Theorem \ref{exactcodim} shows that for $r>2$ zeros along the bad locus are unavoidable whenever the bad locus is nonempty.

\begin{proposition}\label{not volume form}
Let $b\in B$. Suppose that every irreducible component of $B$ through $b$ has codimension at least $3$. Then the bosonic reduction of the extended measure of Theorem \ref{main} vanishes at $b$. In particular, the extended measure is not a volume form near $b$.
\end{proposition}

\begin{proof}
We work on the bosonic reduction near $b$. Let $R=\mc O_{M_{\bos},b}$, which is a regular local ring, and let $L$ denote the universal twisted spin line bundle. Represent $R\pi_*L$ locally by a two-term complex of finite free $R$-modules
\[
F^0\xra{d}F^1.
\]
Set
\[
K:=\ker d=\pi_*L,
\qquad
Q:=\on{coker}d=R^1\pi_*L.
\]
Thus
\begin{equation}\label{K resolution}
0\longrightarrow K\longrightarrow F^0\xra{d}F^1\longrightarrow Q\longrightarrow0.
\end{equation}
By Serre duality,
\[
H^1(C,L)\cong H^0(C,L(-D))^\vee,
\]
so the support of $Q$ is precisely the bad locus $B$.

We first record that $K$ is reflexive. Put $I=\on{im}d$. Then \eqref{K resolution} breaks into
\[
0\to I\to F^1\to Q\to0,
\qquad
0\to K\to F^0\to I\to0.
\]
If $\mathfrak p\notin\on{Supp}Q$, then $K_{\mathfrak p}$ is free. If $\mathfrak p\in\on{Supp}Q$, then $\dim R_{\mathfrak p}=\on{ht}(\mathfrak p)\ge3$. The depth lemma applied to the first sequence gives
\[
\on{depth} I_{\mathfrak p}\ge1,
\]
and applied to the second gives
\[
\on{depth} K_{\mathfrak p}\ge2.
\]
Since $K\subset F^0$, it is torsion-free. Thus $K$ satisfies $S_2$ and is torsion-free over the normal ring $R$, hence is reflexive.

On the good locus, the odd period map is the residue map
\[
\tau:K|_U\longrightarrow \Lambda_1\otimes\mc O_U,
\]
and is injective. We claim that its image cannot have a locally free rank-$r$ extension across $b$. Indeed, if such an extension existed, it would be reflexive and would agree with the reflexive module $K$ away from a subset of codimension at least $2$. Uniqueness of reflexive extension would imply that $K$ itself is locally free, hence free over $R$. But then \eqref{K resolution} would give a free resolution of $Q$ of length at most $2$, so
\[
\on{pd}_R Q\le2.
\]
Auslander--Buchsbaum would imply
\[
\on{depth}_R Q\ge\dim R-2,
\]
and hence $\on{Supp}Q$ would have a component of codimension at most $2$, contradicting the hypothesis.

Now let
\[
\lambda:=\det R\pi_*L.
\]
Over $U$ we have $\lambda|_U=\det K|_U$, so taking the top exterior power of the odd period map gives a Pl\"ucker section
\[
\sigma_U:=\bigwedge^r\tau
\in
\Gamma\left(U,\bigwedge^r(\Lambda_1\otimes\mc O_U)\otimes\lambda^{-1}\right).
\]
Since $B$ has codimension at least $2$, $\sigma_U$ extends uniquely by Hartogs to a section $\sigma$ near $b$. If $\sigma(b)\ne0$, then after shrinking around $b$ it is nowhere zero. The Pl\"ucker relations hold on the dense open set $U$ and therefore everywhere, so $\sigma$ defines a holomorphic map to the Grassmannian $\operatorname{Gr}(r,\Lambda_1\otimes\bb)$ whose pullback of the tautological bundle is a locally free extension of $\on{im}\tau$. This contradicts the preceding paragraph. Hence
\[
\sigma(b)=0.
\]

On the good product, restrict the local formula in the proof of Proposition
\ref{direct extension proposition} to the bosonic reduction.  The mixed blocks
vanish there, so
\[
q_{\bos}=\frac{\det D}{\det A}.
\]
Here $\det A$ is a unit near the conjugate diagonal.  On the same good product,
$\det D$ is the pairing induced by $J_1$ on the two Pl\"ucker sections:
\[
\det D=J_1^{\wedge r}(\sigma_U,\overline{\sigma_U})
\]
in compatible local frames.  The Pl\"ucker sections extend across the bad
locus to $\sigma$ and $\overline\sigma$, so the right-hand side has a
holomorphic extension.  On the other hand,
$\widetilde q_{\bos}\det A$ is a holomorphic extension of the left-hand side
by Proposition \ref{direct extension proposition}.  Since the two extensions
agree on the dense good locus, uniqueness of holomorphic extension gives
\[
\widetilde q_{\bos}\det A
=J_1^{\wedge r}(\sigma,\overline\sigma)
\]
near $(b,\overline b)$.  Because $\sigma(b)=0$ and $\det A$ is a unit,
$\widetilde q_{\bos}(b,\overline b)=0$.  Therefore the bosonic reduction of
\[
\widehat t=\widetilde q^{\,5}e
\]
vanishes at $(b,\overline b)$, and its pullback to the smooth diagonal
vanishes at $b$.
\end{proof}

\begin{remark}\label{conjectural zeros remark}
By Theorem \ref{exactcodim}, every component of the bosonic bad locus $B_{\bos}$ has codimension exactly $r$. Thus, if $B$ is nonempty, then for every $r>2$ Proposition \ref{not volume form} applies along the bad locus: the measure of Theorem \ref{main} extends smoothly, but it necessarily has zeros there and hence is not a volume form.
\end{remark}

\begin{appendices}


\section{Super symplectic pairing} A  super symplectic pairing on a super vector space $V$ is a non-degenerate $\bb$-bilinear form $\langle  \ , \  \rangle: V \times V \to \bb$ such that 
  \begin{itemize}
    \item the restriction $\langle  \ , \rangle_0: V_0 \times V_0 \to \bb$ to the even component $V_0$ is non-degenerate and skew-symmetric, 
    \item the restriction of $\langle \ , \  \rangle_1: V_1 \times V_1 \to \bb$ to the odd component $V_1$ is non-degenerate and symmetric, 
     \item the restriction of $\langle \ , \ \rangle_{1,0}$ and $\langle \ , \ \rangle_{0,1}$ to $V_1 \times V_0$ and $V_0 \times V_1$, respectively,  is identically zero.  \end{itemize}

Let $J:= \langle \ , \ \rangle$ be a supersymplectic pairing on a super vector space $V$. We say that a subspace $W \subset V$ is \emph{Lagrangian} (or, \emph{maximally-isotropic}) with respect to $J$ if the following hold: 
\begin{itemize} 
\item $W_0 \subset V_0$ is Lagrangian (in the usual sense) with respect to the symplectic pairing $J_0$, and 
\item $W_1 \subset V_1$ is maximally-isotropic (in the usual sense) with respect to the symmetric bilinear pairing $J_1$. 
\end{itemize} 

\section{Adding NS punctures}  \label{appendix: ns punctures}

Let $M$ denote the supermoduli space with $2r$ Ramond punctures and $n$ NS punctures, and let
\[
\pi:\mc X\longrightarrow M
\]
be the universal super Riemann surface.  Write $\mc R$ for the universal Ramond divisor.  For the NS punctures, let
\[
s_i:M\longrightarrow \mc X,\qquad i=1,\ldots,n,
\]
be the marked sections, let $P_i\subset\mc X$ denote the associated $0|1$ NS divisors, and set
\[
\mc N:=P_1+\cdots+P_n.
\]
Thus the marked sections and the associated divisors are distinct objects; locally, if an NS puncture is at $z=\theta=0$, the corresponding divisor is $z=0$.

Let $\mc A_X$ denote the sheaf of superconformal vector fields on a fiber $X$.  The tangent space to $M$ at $X$ is
\[
T_XM=H^1\bigl(X,\mc A_X(-N)\bigr).
\]
For a super Riemann surface with Ramond divisor $R$ one has the standard identifications
\[
\mc A_X\cong\mc D^2,
\qquad
\mc D^2\cong\omega_X^{-2}(-2R),
\]
where $\mc D$ is the superconformal distribution.  Hence
\[
\mc A_X(-N)\cong\mc D^2(-N).
\]
Serre duality therefore gives
\[
H^1\bigl(X,\mc A_X(-N)\bigr)^*
\cong
H^0\bigl(X,\omega_X^3(2R+N)\bigr),
\]
and consequently
\begin{equation}\label{canonical-with-NSR}
\omega_M
\cong
\Ber\!\left(R\pi_*\omega_{\mc X/M}^3(2\mc R+\mc N)\right).
\end{equation}

\paragraph{Removing the Ramond polar part.}
The twist by $2\mc R$ must be treated by the sheaf of polar parts, not by ordinary restriction to $\mc R$.  Define
\[
\mc J_{\mc R}(\mc N)
:=
\omega_{\mc X/M}^3(2\mc R+\mc N)
\big/
\omega_{\mc X/M}^3(\mc N).
\]
Then there is a short exact sequence
\[
0\longrightarrow
\omega_{\mc X/M}^3(\mc N)
\longrightarrow
\omega_{\mc X/M}^3(2\mc R+\mc N)
\longrightarrow
\mc J_{\mc R}(\mc N)
\longrightarrow0.
\]
Since $\mc N$ and $\mc R$ are disjoint, the canonical section of
$\mc O(\mc N)$ is nowhere zero in a neighborhood of $\mc R$, and hence gives a canonical identification
\[
\mc J_{\mc R}(\mc N)
\cong
\mc J_{\mc R}(0)
:=
\omega_{\mc X/M}^3(2\mc R)
\big/
\omega_{\mc X/M}^3.
\]
Appendix C of \cite{smiwithRamondwitten2015notes} gives a canonical trivialization
\[
\Ber\!\left(R\pi_*\mc J_{\mc R}(0)\right)\cong\mc O_M.
\]
Multiplicativity of the Berezinian in the resulting exact triangle therefore gives
\begin{equation}\label{finaldescofcotangentbundle}
\omega_M
\cong
\Ber\!\left(R\pi_*\omega_{\mc X/M}^3(\mc N)\right).
\end{equation}

\paragraph{Super Mumford isomorphism with NS and Ramond punctures.}
Consider now the short exact sequence
\[
0\longrightarrow \omega_{\mc X/M}^3
\longrightarrow \omega_{\mc X/M}^3(\mc N)
\longrightarrow \omega_{\mc X/M}^3(\mc N)|_{\mc N}
\longrightarrow0.
\]
Each $P_i$ is finite over $M$, so the restriction term has no higher direct images.  Since the $P_i$ are disjoint, multiplicativity of the Berezinian and \eqref{finaldescofcotangentbundle} give
\begin{equation}\label{stepsmi}
\omega_M
\cong
\Ber\!\left(R\pi_*\omega_{\mc X/M}^3\right)
\otimes
\bigotimes_{i=1}^n \mc L_i,
\end{equation}
where
\[
\mc L_i
:=
\Ber\!\left(
\pi_*\bigl(\omega_{\mc X/M}^3(\mc N)|_{P_i}\bigr)
\right),
\]
so that
\[
\Ber\!\left(
\pi_*\bigl(\omega_{\mc X/M}^3(\mc N)|_{\mc N}\bigr)
\right)
\cong
\bigotimes_{i=1}^n\mc L_i.
\]

Let
\[
f:M\longrightarrow M_0
\]
be the morphism which forgets the NS punctures, where $M_0$ denotes the supermoduli space with the same Ramond punctures and no NS punctures.  In contrast to Ramond punctures, NS punctures are forgettable.  The universal supercurve over $M$ is obtained by base change from the universal supercurve
\[
\pi_0:\mc X_0\longrightarrow M_0,
\]
so
\[
\Ber\!\left(R\pi_*\omega_{\mc X/M}^3\right)
\cong
f^*\Ber\!\left(R\pi_{0*}\omega_{\mc X_0/M_0}^3\right)
\cong
f^*\omega_{M_0}.
\]
Applying the super Mumford isomorphism on $M_0$ now gives
\begin{equation}\label{smiwithNS}
\boxed{
\omega_M
\cong
f^*\!\left(
\Ber^5\!\left(R\pi_{0*}\omega_{\mc X_0/M_0}\right)
\right)
\otimes
\bigotimes_{i=1}^n\mc L_i.}
\end{equation}
The additional factors $\mc L_i$ are the NS insertion-line factors.  Thus this appendix identifies the canonical line on supermoduli with the pullback of the Ramond super Mumford line together with the expected NS puncture factors; producing a scalar measure requires, in addition, the appropriate insertions in these lines.

\section{Even and odd Periods: Explicit description in coordinates}

There are two notions of periods on a super Riemann surface $X$ with Ramond punctures: the periods of closed one-forms on $X$,  and the periods of global sections of the \emph{extended Berezinian sheaf}, $\omega_X'$. 
There is a canonical isomorphism $Z_X^1 \cong H^0(X, \omega_X')$ (cf. \cite{witten2015super}) under which the two notions of periods become equivalent.

\paragraph{Periods of closed one-forms.} The usual $2g$ even periods of a closed one-form $\omega \in Z_X^1$ are given by integrating over a choice of $A$ and $B$ cycles for $H_1(X, \z)=H_1(C,\z)$. 

Its $2r$ odd periods are defined in \cite{witten2015super} as follows: A general closed one form $\omega \in Z_X^1$ is locally in the superconformal coordinates $(z, \theta)$ around a single Ramond puncture of the form \[ \left( f(z) + \frac{\partial g(z)}{\partial z} \theta \right) dz + g(z) d \theta, \]  where  $g$ and $f$ are local holomorphic functions on $X$, possibly depending on parameters from the base.  
 Let $j: R \to X$ be the natural inclusion of the Ramond divisor, and let $x:=j^*(\theta)$. 
The restriction of $\omega$ to the Ramond puncture contained in $(z,\theta)$ is $g(0) dx $. Witten defines the constant $g(0)$ to be the odd period of $\omega$. Note that $g(0) dx$ is a closed one-form on $R$. Repeating this for every component of $R$, we get the $2r$-tuple of odd period of $\omega$: 
 \begin{equation} \omega \vert_R =(g_1(0) dx^1, \dots, g_{2r}(0) dx^{2r}) \in Z_R^1 \end{equation}
The constants are unique up to sign. 

More abstractly, the map computing the odd periods of closed one-forms on $R$ is given by the natural restriction map $Z_X^1 \to Z_R^1$ of closed one-forms on $X$ to closed one-forms on $R$. 

\paragraph{Periods of global sections of the extended Berezinian sheaf.}

We define the periods of sections of $H^0(X, \omega')$ in Section \ref{sectiontitle: super period map}. Let us give an explicit description of their odd periods in terms of local superconformal coordinates $(z, \theta)$ near a Ramond puncture. 

A general global section $\omega \in H^0(X, \omega')$ is locally in the coordinates $(z, \theta)$ given by 
\[ \omega= \left( \frac{g(z)}{z} + f(z) \theta \right)[dz| d \theta]. \]
Restricting $\omega$ to $R$ is equivalent to computing the residue of the function $\left( \frac{g(z)}{z} + f(z) \theta \right)$ at $z=0$, and hence $\omega \vert_R=g(0) dx$. We define $g(0)$ to be the odd period of $\omega$ at the Ramond puncture contained in $(z, \theta)$.  Applying this to every Ramond puncture,  we get the $2r$-tuple of odd periods of $\omega$. 

The fact that these two notions of periods are equivalent  is immediate from the local description of the canonical isomorphism $Z_X^1 \cong H^0(X, \omega')$ in the coordinates $(z, \theta)$: 
\[   \left( f(z) + \frac{\partial g(z)}{\partial z} \theta \right) dz + g(z) d \theta \mapsto \left( \frac{g(z)}{z} + f(z) \theta \right) [dz| d \theta] \]

\end{appendices}

\providecommand{\bysame}{\leavevmode\hbox to3em{\hrulefill}\thinspace}
\providecommand{\MR}{\relax\ifhmode\unskip\space\fi MR }
\providecommand{\MRhref}[2]{%
  \href{http://www.ams.org/mathscinet-getitem?mr=#1}{#2}
}
\providecommand{\href}[2]{#2}

\end{document}